\documentclass[11pt]{article}

\usepackage{latexsym}

\usepackage{amsmath}
\usepackage{amssymb}

\def\be{\begin{equation}}
\def\ee{\end{equation}}
\def\ba{\begin{eqnarray*}}
\def\ea{\end{eqnarray*}}
\def\baa{\begin{eqnarray}}
\def\eaa{\end{eqnarray}}

\def\e{\varepsilon}
\def\a{\alpha}
\def\la{\lambda}

\title{\Large\bf A note on a new cubically convergent one-parameter root solver}

\author{{\bf L. D. Petkovi\'c$^{\,1,}$}\footnote{Corresponding author}\;,\;\bf{M. S. Petkovi\'c$^{\,2,}$}
   \\[2mm]
$^1$\small\it   Faculty of Mechanical Engineering,
University of Ni\v s,\\[-1mm]
\small\it A. Medvedeva 14, 18000 Ni\v s, Serbia\\
$^2$\small\it  Faculty of Electronic Engineering,
University of Ni\v s, \\[-1mm] \small \it A. Medvedeva 14, 18000 Ni\v s, Serbia}

\date{}

\begin{document}

\maketitle

\begin{abstract}
A new one-parameter family of iterative method for solving nonlinear equations is constructed and studied. Two variants, both with cubic convergence, are developed, one for finding  simple zeros and other for  multiple zeros of known multiplicities. This family generates a variety of different third order methods, including Halley-like method as a special case. Four numerical examples are given to demonstrate convergence properties of the  proposed methods  for multiple zeros and various values of the parameter.
\\[1mm]
\end{abstract}

\renewcommand{\thefootnote}{}
\footnote{{\it E-mail addresses:} ljiljana@masfak.ni.ac.rs (L. D. Petkovi\'c), msp@junis.ni.ac.rs (M. S. Petkovi\'c)
 }

\section{Introduction}

Approximating zeros of a given scalar function $f$ belongs to the most important problems that occur not only in applied mathematics but also in many disciplines of engineering branches, computer science, physics, finance,  and so on. 
Since there is a vast  number of papers and books devoted to iterative methods for finding simple and multiple roots of nonlinear equations, see, e.g., \cite{traub}--\cite{neta2},  we will not discuss in details characteristics of existing methods.

The main goal of this paper is to present two new one-parameter families of iterative methods for finding simple or multiple zeros of a given function. The main advantages of this family are: 1) the ability to generate a variety of different cubically convergent methods;  the proposed family can serve for the construction of very fast iterative methods for approximating all zeros of a polynomial, see \cite{AML}.

The paper is organized as follows.
In Section 2 we construct a one-parameter family of iterative methods for finding simple roots of nonlinear equations and show that its order of convergence is three.  In Section 3 the iterative formula for simple zeros is directly used for the construction of a one-parameter family for finding multiple zeros of the known multiplicity. Results of numerical experiments for several values of the parameter $p$ through three iteration steps are displayed in Section 4 using four test functions.

\section{One-parameter family for simple zeros}

We begin this section with Traub's result given in \cite[Theorem 2.5]{traub}.

\bigskip

{\bf Theorem 1.} {\it Let $\psi(x)$ be an iteration function which defines an iterative method for finding a zero $\a$ of multiplicity $m$ of a given function $f.$   Then for these values of $m$ there exists a function $\omega(x)$ such that
  \be
  \psi(x)=x-u(x)\omega(x),\quad u(x)=\frac{f(x)}{f'(x)},\quad \omega(\a)\ne 0.\label{1}
  \ee}

  In this paper we will restrict our attention to iterative methods with cubic convergence ($r=3$). We will often use the abbreviations
  $$
  A_{\la}(x)=\frac{f^{(\la)}(\a)}{\la!f'(\a)}\quad (\la=1,2,\ldots).
  $$ For brevity,  we will write sometimes only $u$ instead of $u(x)$ for short. The abbreviation $AEC({\rm IM})$ will denote {\it asymptotic error constant} of the iterative method (IM). First we present two well known cubically convergent  methods  free of squares:
  \baa
\hspace*{-0.9cm}&&C(x)=x-u(x)\Bigl(1+A_2(x)u(x)\Bigr)\ \ \mbox{\rm (Chebyshev's method)},\label{2}\\
\hspace*{-0.9cm}&&  H(x)=x-\dfrac{u(x)}{1-A_2(x)u(x)}\ \  \mbox{\rm (Halley's method)}.\label{3}
\eaa

Regarding  (\ref{1})  we note that $\omega(u)=1+A_2u$ for Chebyshev's method (\ref{2}) and $\omega(u)=1/(1-A_2u)$ for Halley's method (\ref{3}). Therefore, $\omega(u)$ is a polynomial approximation  in (\ref{2}), while $\omega(u)$ is a rational  approximation in (\ref{3}).
In this paper we will consider a rational  approximation to construct a new cubically convergent iterative method in the form
 \be G(u)=x-u(x)\cdot \frac{a+p\cdot u(x)}{1+c\cdot u(x)}.\label{7}
  \ee
We allow that the coefficients $a$ and $c$  in (\ref{7}) can be constants as well as  some functions of the argument $x,$ while $p$ is a real or complex parameter.

First, we start from  Schr\"oder-Traub {\it basic sequence} $\{E_k\}$ defined recursively by
\be
\left\{\begin{array}{l}
E_2(x)=x-u(x),\\[8pt]
E_{k+1}(x)=E_k(x)-\dfrac{u(x)}{k}E'_k(x)\quad (k\ge 2),\end{array} \right. \label{8}
\ee
which defines the generalized iterative method of order $k+1$ in the form of a power series, see Traub \cite[Sec. 5.1]{traub}. For example (suppressing the argument $x$)
\ba
E_3(x)&=&x-u-A_2u^2\quad \mbox{\rm (Chebyshev's method (\ref{2})}),\\
E_4(x)&=&x-u-A_2u^2-(2A_2^2-A_3)u^3,\\
E_5(x)&=&x-u-A_2u^2-(2A_2^2-A_3)u^3-(5A_2^3-5A_2A_3+A_4)u^4,\ \ \mbox{\rm etc.}
\ea
We will employ the following assertion.

\bigskip

{\bf Theorem 2.} (Schr\"oder \cite{sreder}) {\it Any root-finding algorithm $F_n$ of the order $n$ can be presented in the form
 \be
 F_n(x)=E_n(x)+f(x)^n \eta_n(x),\label{9}
 \ee
 where $\eta_n$ is a function bounded in $\a$ and depending on $f$ and its derivatives.}

 \medskip

Let $G_3(x;p)$ be the root-solver  to be constructed. According to Theorem 2 we seek for the coefficients $a,\ p,\ c$ in (\ref{7})
so that
 \be
 G_3(x;p)=E_3(x)+f(x)^3 \eta_3(x)\label{10}
 \ee
 holds.

 For two real or complex numbers $z$ and $w$ we will write $z=O_M(w)$ if $|z|=O(|w|)$ (the same order of their moduli),
where $O$ represents the Landau symbol.
After the development in geometric series we have
 $$
 G_3(u)=x-u\cdot \frac{a+pu}{1+cu}=x-u(a+pu)(1-cu+c^2u^2+\cdots)=x-au+(ac-p)u^2+\cdots\ .
 $$
 Using this relation and (\ref{10}), and applying the method of undetermined coefficients, we obtain $a=1$ and  $c=p-A_2.$ In this way we have constructed the following one-parameter cubically convergent iterative method
 \be
 \hat x=G_3(x;p)=x-\frac{u(x)\Bigl(1+p\,u(x)\Bigr)}{1+\Bigl(p-A_2(x)\Bigr)u(x)},\label{11}
 \ee
 where $\hat x$ is a new approximation to the zero $\a$ of $f.$

Let $\e=x-\a$ be the approximation error. To find asymptotic error constant ($AEC$ for short), we use the developments in Taylor' series:
\ba
f(x)&=&f'(\a)\Bigl(\e+A_2\e^2+A_3\e^3+O_M(\e^4)\Bigr),\\
f'(x)&=&f'(\a)\Bigl(1+2A_2\e+3A_3\e^2+O_M(\e^3)\Bigr),\\
f''(x)&=&f'(\a)\Bigl(2A_2+6A_3\e+O_M(\e^2)\Bigr).
\ea
Hence
\ba
u(x)&=&\frac{f(x)}{f'(x)}=\e-A_2\e^2+(2A_2^2-2A_3)\e^3+O_M(\e^4),\\
A_2(x)&=&\frac{f''(x)}{2f'(x)}=A_2+(3A_3-2A_2^2)\e+(4A_2^3-9A_2A_3)\e^2+O_M(\e^3).
\ea
Substituting the  expressions for $u(x)$ and $A_2(x)$ in (\ref{11}), we obtain
 \be
 \hat \e=\hat x-\a=(A_2^2-A_3+pA_2)\e^3+O_M(\e^4).\label{12}
 \ee

 From (\ref{12}) we immediately state the following assertion.

\bigskip

{\bf Theorem 3.} {\it  Assume that $x_0$ is sufficiently close initial approximation to the zero $\a$ of at least three-time differentiable function $f.$ Then the one-parameter family of iterative methods
 \be
 x_{k+1}=x_k-\frac{u(x_k)\Bigl(1+p\,u(x_k)\Bigr)}{1+\Bigl(p-A_2(x_k)\Bigr)u(x_k)}\quad (k=0,1,\ldots)\label{13}
 \ee
 has the order of convergence three for any real or complex parameter $p,$ bounded in magnitude, and
 $$
 AEC(\ref{13})=\lim_{k\to \infty} \left| \frac{x_{k+1}-\a}{(x_k-\a)^3}\right|=\left|A_2^2(\a)-A_3(\a)+
 pA_2(\a)\right|
  $$
  is valid.
 }

\bigskip
{\bf Remark 1.} In a special case when $p=0$, from (\ref{12}) we obtain Halley's method (\ref{3}) with
 $AEC(\ref{3})=|A_2^2(\a)-A_3(\a)|,$ which is well-known result. Furthermore, when $p\to \pm \infty,$ then the method (\ref{13}) reduce to quadratically convergent Newton's method $x_{k+1}=x_k-u(x_k).$ For this reason, {\it one should avoid the choice of the parameter $p$ large in magnitude}.

\bigskip

{\bf Remark 2.} We restrict ourselves that the parameter $p$ is a real
or complex constant. It is interesting to consider another special case $p=A_2(x_k)$ which could happen accidentally  in the $k$-iteration. Then the iterative process (\ref{13}) switches to Chebyshev's method
$$
 x_{k+1}=x_k-u(x_k)\Bigl(1+A_2(x_k)\,u(x_k)\Bigr)
\quad (k=0,1,\ldots),
$$
see (\ref{2}).
Since the probability of this case is 0, we will not discuss Chebyshev's method in what follows. However, the described case can be helpful in finding suitable range od the parameter $p.$

\bigskip

{\bf Remark 3.} If we choose $p=\bigl(A_3(x_k)-A_2(x_k)^2\bigr)/A_2(x_k)$ in (\ref{13}),  the iterative method (11) becomes
\be
x_{k+1}=x_k-u(x_k)\biggl(1+\frac{A_2(x_k)u(x_k)}{A_2(x_k)+
(A_3(x_k)-2A_2^2(x_k))u(x_k)}\biggl)\quad (k=0,1,\ldots).\label{red4}
\ee
From (\ref{12}) we find that the iterative method (\ref{red4}) has the order of convergence equal to 4. Developing the denominator of the expression in the above parenthesis
   around the point $u=0,$ the iterative formula (\ref{red4}) reduced to the fourth order method $E_4(x_k)$ given above.

\section{Multiple zeros}

 Let us now consider the case when $\a$ is the zero of $f$ of the known order of multiplicity $m\ge 1.$ Note that $\a$ is a simple zero for the function
  $$
  F(x)=f(x)^{1/m}.
  $$
     We find the first two derivatives of $F:$
   \be
   F'(x)=\frac{F(x)f'(x)}{mf(x)},\quad F''(x)=F(x)\cdot \frac{(1-m)f'(x)^2+mf(x)f''(x)}{m^2f(x)^2}.\label{14}
   \ee

  Taking into account the expressions for the derivatives $F'$ and $F''$ given by (\ref{14}), let us  replace  $u(x)$ and $A_2(x),$ appearing in (\ref{11}),  by new functions $v(x)$ and $d_2(x)$ defined by
   \be
   v(x):=\frac{F(x)}{F'(x)}=\frac{mf(x)}{f'(x)},\quad d_2(x):=\frac{F''(x)}{2F'(x)}=\frac{(1-m)f'(x)^2+mf(x)f''(x)}
   {2mf(x)f'(x)}.
   \label{15}
   \ee
Then the iteration function (\ref{11}) becomes
\be
\hat x=G_m(x;p)=x-\frac{v(x)\bigl(1+pv(x)\bigr)}{1+\bigl(p-d_2(x)
\bigr)v(x)}.\label{16}
 \ee

Replacing the expressions (\ref{15}) for $v(x)$ and $d_2(x)$  we can modify the iteration function (\ref{11}) for simple zeros to the following iteration function for finding multiple zeros
\be
x_{k+1}=G_m(x_k;p)=x_k-\frac{2mu(x_k)\bigl(1+mpu(x_k)\bigr)}
{1+m+2m\bigl(p-A_2(x_k)
\bigr)u(x_k)}\quad (k=0,1,\ldots).\label{17}
 \ee
 Note that the choice of $p=0$ in (\ref{17}) gives Halley-like method for finding multiple zeros \cite{bodevig}
\be
x_{k+1}=x_k-\frac{u(x_k)}{\dfrac{m+1}{2m}-A_2(x_k)u(x_k)}\quad (k=0,1,\ldots).\label{halej}
\ee

 \bigskip

\medskip
{\bf Theorem 4.} {\it Let $x_0$ be sufficiently close initial approximation to the zero $\a$ of the known multiplicity $m\ge 1$ of a given function $f.$ Then the iterative method {\rm (\ref{17})} is cubically convergent and
\be
AEC(\ref{17})=\lim_{k\to \infty} \left|\frac{x_{k+1}-\a}{(x_k-\a)^3}\right|=\left|\frac{pB_{m+1}}{mB_m}-\frac{B_{m+2}}{mB_m}
+\frac{(m+1)B_{m+1}^2}{2m^2B_m^2}\right|\label{18}
\ee
is valid,
where $B_r=f^{(r)}(\a)/r!\; .$
}

\medskip

{\bf Proof.} To prove the theorem we will use the iteration function (\ref{16}). Introduce the errors of approximations $\e=x-\a,\ \hat \e=\hat x-\a$ and coefficients
 $$
 C_r=\frac{m!}{(m+r)!} \frac{f^{(m+r)}(\a)}{f^{(m)}(\a)}\quad (r=1,2,\ldots).
 $$
 Then the following developments in Taylor serious are valid:
 \ba
 f(x)&=&B_m\e^m\Bigl(1+C_1\e+C_2\e^2+C_3\e^3+O_M(\e^4)\Bigr),\\
 f'(x)&=&B_m\e^{m-1}\Bigl(m+(m+1)C_1\e+(m+2)C_2\e^2+(m+3)C_3\e^3+O_M(\e^4)\Bigr),\\
 f''(x)&=&B_m\e^{m-2}\Bigl(m(m-1)+m(m+1)C_1\e+(m+1)(m+2)C_2\e^2\\
 && \hskip1.3cm+(m+2)(m+3)C_3\e^3+O_M(\e^4)\Bigr).
 \ea
 Using (\ref{15}) and these expressions, we find
 \be
 \left\{\begin{array}{l} v(x)=\e-\dfrac{C_1\e^2}{m}+\dfrac{((m+1)C_1^2-2mC_2)\e^3}{m^2}+O_M(\e^4),
\\[8pt]
 d_2(x)=\dfrac{C_1}{m}-\dfrac{((m+1)C_1^2-6mC_2)\e}{2m^2}+O_M(\e^2).
 \end{array}\right.\label{19}
 \ee

 According to (\ref{19}) and Taylor's series it follows
 $$
 \frac{1}{1+(p-d_2(x))v(x)}=
  1+\left(\frac{{C_1}}{m}-p\right)\e+\frac{
   \left(- {C_1}^2(3m+1)-2mp {C_1}
   +6m {C_2}+2 m^2 p^2\right)\e^2}{2
   m^2}+O\left(\e^3\right).
   $$
   Using the last expression and (\ref{16}), we find after short arrangement
   \be
   \hat \e=\hat x-\a=\frac{\Bigl((m+1)C_1^2+2mpC_1-2mC_2\Bigr)\e^3}{2m^2}+O_M(\e^4).\label{20}
   \ee

   With regard to (\ref{20}) it follows that the order  of the iterative method (\ref{17}) is three. Since $C_r=B_{m+r}/B_m,$ from (\ref{20}) we obtain
   the asymptotic error constant $AEC({\rm \ref{17})},$ given by (\ref{18}).\ $\square$
   \bigskip

   \section{Numerical results}

    The theoretical order of convergence  of the iterative method (\ref{17}) is  three, see Theorem 4. However, it is always convenient to check the convergence behavior in practice. For this reason, in our numerical experiments we have calculated the so-called {\it computational order of convergence $r_c$} (COC, for brevity) using  the approximate formula
 \be
 r_c=\frac{\log|f(x_{k+1})/f(x_k)|}{\log|f(x_{k})/f(x_{k-1})|}.\label{COC}
 \ee
Note that the formula (\ref{COC}) is a special case of a general formula given in \cite{jay}. The tested functions are given in Table \ref{t1}.

 \begin{table}[htb]
 \begin{center}
  {\small
    \begin{tabular}{|lc|l|c|l|}
 \hline
  $f(x)$ & &$m$  & $x_0$ & $\a$\\ \hline
  $f_1(x)=\bigl(x\,\sin x-2\sin^2(x/\sqrt{2})\bigr)\bigl(x^5+x^2+100\bigr)$ & & $6$  & $-1.2$ & $0$\\[4pt]
 $f_2(x)=\bigl(xe^{x^2}-\sin^2x+3\cos x+5)^2$ & & $2$  & $-1$ & $-1.2076478271309\ldots$\\[4pt]
 $f_3(x)=\bigl(e^{x^2+4x+5}-1\bigr)^3\sin^2(t+2-i)$ & & $5$  & $-1.7+0.8i$ & $-2+i$\\[4pt]
  $f_4(x)=\bigl(x-\sin x\bigr)^4$  & & $12$ & $0.4$ & $0$\\ \hline
    \end{tabular}
 \caption{Tested functions for $f_1-f_4$\label{t1}}
 }
 \end{center}
 \end{table}

\begin{table}[htb]
\begin{center}
{
\begin{tabular}{|c|l|l|l|l|}\hline
\multicolumn{5}{|c|}{ $f_1(x)=\bigl(x\,\sin x-2\sin^2(x/\sqrt{2})\bigr)\bigl(x^5+x^2+100\bigr)$}\\ \hline\hline
$p$ & $|x_1-\a|$ &  $|x_2-\a|$ & $|x_3-\a|$ & $r_c$ (\ref{COC})\\ \hline
$-2$ & $2.29(-2)$ & $1.40(-7)$ & $2.84(-23)$ & $3.011$ \\ \hline
$-1$ & $8.91(-4)  $ & $7.25(-12)$ & $\setlength{\fboxrule}{1.1pt}\fbox{$3.90(-36)$} $ & $ $3.000$  $ \\ \hline
$0$ & $7.08(-2) $ & $3.64(-6) $ & $ 3.39(-19) $ & $ 3.000  $ \\ \hline
$1$ & $ 0.111 $ & $ 1.42(-2) $ & $ 3.06(-8) $ & $ 3.000  $ \\ \hline
$2$ & $ 0.172 $ & $1.19(-5)  $ & $ 1.72(-17) $ & $2.846  $ \\ \hline\hline
\multicolumn{5}{|c|}{$f_2(x)=(xe^{x^2}-\sin^2x+3\cos x+5)^2$}\\ \hline\hline
$-2$ & $4.93(-2)$ & $4.34(-4)$ & $2.66(-10)$ & $3.067$ \\ \hline
$-1$ & $1.87(-2)  $ & $1.17(-5)  $ & $ 2.82(-15) $ & $ $3.013$  $ \\ \hline
$0$ & $7.99(-4)  $ & $1.29(-10)  $ & $ \setlength{\fboxrule}{1.1pt}\fbox{$5.50(-31)$} $ & $ 3.000  $ \\ \hline
$1$ & $ 1.10(-2) $ & $ 1.65(-6) $ & $ 5.64(-18) $ & $ 2.994  $ \\ \hline
$2$ & $ 1.93(-2) $ & $2.04(-5)  $ & $ 2.32(-14) $ & $2.991   $ \\ \hline\hline
\multicolumn{5}{|c|}{$f_3(x)=(e^{x^2+4x+5}-1)^3\sin^2(t+2-i)$}\\ \hline\hline
$-2$ & $6.17(-2)$ & $1.74(-4)$ & $3.45(-12)$ & $3.031$ \\ \hline
$-1$ & $3.30(-2) $ & $1.44(-5)  $ & $ 1.18(-15) $ & $ $3.007$  $ \\ \hline
$0$ & $1.33(-2)  $ & $2.94(-7)  $ & $ 5.32(-20) $ & $ 3.000  $ \\ \hline
$1$ & $ 7.04(-2) $ & $ 1.36(-7) $ & $ \setlength{\fboxrule}{1.1pt}\fbox{$9.83(-22)$} $ & $ 2.999  $ \\ \hline
$2$ & $ 1.06(-2) $ & $7.59(-7)  $ & $ 2.85(-19) $ & $2.997   $ \\ \hline\hline
\multicolumn{5}{|c|}{$f_4(x)=\bigl(x-\sin x\bigr)^4$}\\ \hline\hline
$-2$ & $1.38(-2)$ & $4.47(-8)$ & $1.78(-24)$ & $3.067$ \\ \hline
$-1$ & $3.21(-3)  $ & $5.59(-10)  $ & $ 2.91(-30) $ & $ $3.001$  $ \\ \hline
$0$ & $1.08(-3)  $ & $2.08(-11)  $ & $ 1.50(-34) $ & $ 3.000  $ \\ \hline
$1$ & $ 1.58(-4) $ & $ 6.52(-14) $ & $ \setlength{\fboxrule}{1.1pt}\fbox{$4.63(-42)$} $ & $ 3.000  $ \\ \hline
$2$ & $ 3.53(-4) $ & $7.37(-13)  $ & $ 6.68(-39) $ & $3.000   $ \\ \hline
\end{tabular}
\caption{Errors of approximations; functions $f_1-f_4$\label{t2}}
}
\end{center}
\end{table}

In Table \ref{t2}  we have presented the errors of approximations $\e_k=|z_k-\a|\ (k=1,2,3)$ produced by the  method (\ref{17}) for 5 values of the parameter $p.$ The denotation $A(-h)$ means $A\times 10^{-h}.$
The most accurate approximations, obtained after the third iterative step, are boxed  in Table \ref{t2}. We observe that the best results are obtained taking $p=-1$ for $f_1$, $p=0$ for $f_2,$  and $p=1$ for $f_3$ and $f_4.$

Except the functions listed in Table \ref{t1}, we have also tested a number of functions of various structure.  However, we have not found the value of $p$ which defines approximately the best method from the family (\ref{17}). The influence of the parameter $p$ to the accuracy of approximations to the zeros of a given function is very complex and it is hard to find its optimal value even within a
particular class of functions. From the  discussion given in Remark 1 we can conclude that large values of $p$ are not convenient since the order of convergence decreases and tends to 2. Furthermore, for $p=0$ the method (\ref{13}) reduces to Halley's method which belongs to the group of cubically convergent methods with very good convergence behavior.

Our numerical experiments have shown that optimal parameter $p$ for some classes of functions takes negative values belonging to the interval $[-b,0]\ (b>0).$ According to all facts mentioned above and
a number of tested functions (see Table \ref{t2} for demonstration), we have concluded  that $p$ should be taken from the interval $[-a,a]\ (a>0)$ for relatively small $a,$ say, $a\le 3.$ Following Remark 2, there follows that the choice of $p$ very close to $A_2(x_k)$ could be also good (taking $p\approx A_2(x_k)$ before running the $k$-iteration). However, if an initial approximation $x_0$ is not sufficiently close to the sought zero, the values $A_2(x_0)$ can be rather crude, producing slow convergence at the beginning of iterative process.

\bigskip

{\bf Remark 4.} The values of COC $r_c$ in Table \ref{t2}, taken with 3 decimal digits of mantissa (thus, 3.000 does not mean 3) mainly match well the theoretical order 3. However, in some cases unexpected values of $r_c$ appear.  The explanation is simple: formula (\ref{COC}) works well when the approximations $x_{k-1},\;x_k,\;x_{k+1}$ are sufficiently close to the zero. One additional iteration more would give more realistic value of $r_c.$

\end{document}